\documentclass[12pt]{article}
\usepackage{amssymb}

\newtheorem{thm}     {Theorem}[section]

\newtheorem{cor}     [thm]{Corollary}
\newtheorem{lemma}   [thm]{Lemma}

\newcommand{\proof} {\noindent{\bf Proof. }}

\newcommand{\B}{\mathbb B}
\newcommand{\C}{\mathbb C}
\newcommand{\D}{\mathbb D}

\newcommand{\R}{\mathbb R}

\newcommand{\st}{{\rm st}}

\def\Re{{\rm Re\,}}

\def\bar{\overline}

\begin{document}

\title{Discs and boundary uniqueness for psh functions on an almost complex manifold  }
\author{Alexandre Sukhov{*} }
\date{}
\maketitle

{\small

* Universit\'e des Sciences et Technologies de Lille, Laboratoire
Paul Painlev\'e,
U.F.R. de
Math\'e-matique, 59655 Villeneuve d'Ascq, Cedex, France, sukhov@math.univ-lille1.fr
The author is partially suported by Labex CEMPI.

Instiut of Mathematics with Computing Centre - Subdivision of the Ufa Research Centre of Russian
Academy of Sciences, 45008, Chernyshevsky Str. 112, Ufa, Russia.

}
\bigskip

{\small Abstract. We prove that a totally real submanifold (of maximal dimension) of an almost complex manifold is a boundary uniqueness set for plurisubharmonic functions}.



\bigskip

\section{Introduction}
This paper addresses some aspects of pluripotential theory on almost complex manifolds. In the general case of non-integrable almost complex structures, the development  of this theory began recently and many natural questions remain open. Our main motivation arises from the  paper by J.-P.Rosay \cite{Ro} where several interesting properties of (non) pluripolar subsets of almost complex manifolds are established. In particular, he proved that a $J$-complex curve is a pluripolar set. On the other hand, it was sketched in \cite{Ro} a proof of the fact that a totally real submanifold of dimension 2 can non be a pluripolar subset of an almost complex manifold of complex dimension 2 (a hint for higher dimension argument also was indicated without details).

In the present paper we prove a considerably more general result (Theorem \ref{theo1}) which can be viewed as a boundary uniqueness theorem for plurisubharmonic functions. It states that a function plurisubharmonic in a wedge with generic totally real edge is equal to $-\infty$ identically if it tends to 
$-\infty$ approaching the edge. Our proof is completely different from the  argument of \cite{Ro}. Our approach  is based on construction of a suitable family of $J$-complex discs and is inspired by the well-known work \cite{Pi}. We hope that the almost complex version of this construction presented here will have other applications.  Another approach to filling a wedge by discs is used in \cite{Su} and is based on stability results for the Riemann-Hilbert boundary value problem. The method of present paper is more direct and gives an additional information on geometry of discs. The required family of complex discs is obtained as a solution of a suitable integral equation generalizing the classical Bishop method. In the almost complex case this equation arises from the Cauchy-Green type formula.

The paper is organized as following. In Section 2 we recall basic properties of almost complex manifolds (see more in \cite{Aud}). Section 3 contains the construction of Bishop type discs glued to a totally real manifold along a boundary  arc. In section 4 we deduce the main result on boundary uniqueness for psh functions.

\section{ Almost complex manifolds, discs and plurisubharmonic functions} Everywhere through this paper we assume that manifolds and structures are $C^\infty$ smooth although the main results remain true under considerably weaker regularity assumptions.

Let $M$ be a smooth  manifold of dimension $2n$. {\it An almost complex structure} $J$ on $M$ is a smooth map  which associates to every point $p \in M$ a linear isomorphism $J(p): T_pM \to T_pM$ of the tangent space $T_pM$ satisfying $J(p)^2 = -I$; here  $I$ denotes the identity map of $T_pM$. Thus, $J(p)$ is a linear complex structure on $T_pM$. A couple $(M,J)$ is called {\it an almost complex manifold} of complex dimension n. Note that every almost complex manifolds admits the canonical orientation.

The {\it standard complex structure} $J_{st} = J_{st}^{(2)}$ on $M = \R^2$  is given by the matrix 

\begin{eqnarray}
\label{J_st}
J_{st}= \left(
\begin{array}{cll}
0 & & -1\\
1 & & 0
\end{array}
\right)
\end{eqnarray}
(in the canonical coordinates of $\R^2$. More generally, the standard complex structure $J_{st}$ on $\R^{2n}$ is represented by the block diagonal matrix $diag(J_{st}^{(2)},...,J_{st}^{(2)})$ (we skip the notation of dimension because it will be clear from the context). As usual,  setting $iv := Jv$ for $v \in \R^{2n}$, we identify $(\R^{2n},J_{st})$ with $\C^n$ having 
$z = x + iy = x + Jy$ for the standard complex coordinates $z = (z_1,...,z_n) \in \C^n$.

Let $(M,J)$ and $(M',J')$ be smooth  almost complex manifolds. A $C^1$-map $f:M' \to M$ is called  
$(J',J)$-complex or  $(J',J)$-holomorphic  if it satisfies {\it the Cauchy-Riemann equations} 
\begin{eqnarray}
\label{CRglobal}
df \circ J' = J \circ df.
\end{eqnarray}

Of course,   a map $f: \C^n \to \C^m$ is $(J_{st},J_{st})$-holomorphic if and only if each component of $f$ is a usual holomorphic function.

 Every almost complex manifold
$(M,J)$ can be viewed locally as the unit ball $\B$ in
$\C^n$ equipped with a small almost complex
deformation of $J_{st}$. Indeed, we have the following useful statement.
\begin{lemma}
\label{lemma1}
Let $(M,J)$ be an almost complex manifold. Then for every point $p \in
M$, every real $\alpha \geq 0$ and   $\lambda_0 > 0$ there exist a neighborhood $U$ of $p$ and a
coordinate diffeomorphism $z: U \rightarrow \B$ such that
$z(p) = 0$, $dz(p) \circ J(p) \circ dz^{-1}(0) = J_{st}$  and the
direct image $ z_*(J): = dz \circ J \circ dz^{-1}$ satisfies $\vert\vert z_*(J) - J_{st}
\vert\vert_{C^\alpha(\bar {\B})} \leq \lambda_0$.
\end{lemma}
\proof There exists a diffeomorphism $z$ from a neighborhood $U'$ of
$p \in M$ onto $\B$ satisfying $z(p) = 0$ and $dz(p) \circ J(p)
\circ dz^{-1}(0) = J_{st}$. For $\lambda > 0$ consider the dilation
$d_{\lambda}: t \mapsto \lambda^{-1}t$ in $\R^{2n}$ and the composition
$z_{\lambda} = d_{\lambda} \circ z$. Then $\lim_{\lambda \rightarrow
0} \vert\vert (z_{\lambda})_{*}(J) - J_{st} \vert\vert_{C^\alpha(\bar
{\B})} = 0$ for every real $\alpha \geq 0$. Setting $U = z^{-1}_{\lambda}(\B)$ for
$\lambda > 0$ small enough, we obtain the desired statement.

\bigskip

 Gromov's theory of pseudoholomorphic curves studies the property of solutions  $f$ of (\ref{CRglobal}) in the special case   where $M'$ has the complex dimension 1. These holomorphic maps are called $J$-complex (or $J$-holomorphic) curves.
We use the notation  $\D = \{ \zeta \in \C: \vert \zeta \vert < 1 \}$ for  the
unit disc in $\C$ always assuming that it is equipped with the standard complex structure   $J_{\st}$. If in the above definition we have $M' = \D$  we  call such a map $f$ a $J$-{\it complex  disc} or a
 pseudo-holomorphic disc or just a  holomorphic disc
if $J$ is fixed. Similarly, if $M'$ is the Riemann sphere, $f$ is called a $J$-{\it complex sphere}. These two classes of pseudo-holomorphic curves are particularly useful. 

Let $(M,J)$ be an almost complex manifold and $E \subset M$ be a real submanifold of $M$. 
Suppose that a $J$-complex disc $f:\D \to M$ is  continuous on $\overline\D$.  We some abuse of terminology, we also call the image $f(\D)$  simply by a disc and twe call he image $f(b\D)$ by  the boundary of a disc. If  $f(b\D) \subset E$, then we say that (the boundary of ) the disc  $f$ is {\it glued} or {\it attached} to $E$ or simply 
that $f$ is attached to $E$. Sometimes such maps are called {\it Bishop discs} for $E$ and we employ this terminology. Of course, if $p$ is a point of $E$, 
then the constant map $f \equiv p$ always satisfies this definition. 

In this paper we deal with a special class of real submanifolds. A submanifold $E$ of an almost complex $n$-dimensional $(M,J)$ is called 
{\it totally real} if at every point $p \in E$ the tangent space $T_pE$  does not contain non-trivial complex vectors that is $T_pE \cap JT_pE = \{ 0 \}$. This is well-known that the (real)  dimension of a totally real submanifold of $M$ is not bigger than $n$; we will consider in this paper only $n$-dimensional totally real submanifolds that is the case of maximal dimension.

A totally real manifold $E$ can be defined as 
 \begin{eqnarray}
 \label{edge}
 E = \{ p \in M: \rho_j(p) = 0 \}
 \end{eqnarray}
 where $\rho_j:M   \to  \R$  are smooth functions with non-vanishing gradients.  The condition of total reality means that for every $p \in E$ the $J$-complex linear parts of the differentials $d\rho_j$ are (complex) linearly independent.

 A subdomain

\begin{eqnarray}
\label{wedge}
 W = \{ p \in M: \rho_j < 0, j = 1,...,n \}.
 \end{eqnarray}t
 is called {\it the wedge with the edge} $E$.

\bigskip

\subsection{ Cauchy-Riemann equations in coordinates} Everything will be local, so (as above) we are in a neighborhood $\Omega$ of $0$ in $\C^n$ with the standard complex coordinates $z = (z_1,...,z_n)$. We assume that $J$ is an almost complex structure defined on $\Omega$ and $J(0) = J_{st}$. Let 
$$z:\D \to \Omega,$$ 
$$z : \zeta \mapsto z(\zeta)$$ 
be a $J$-complex disc. Setting $\zeta = \xi + i\eta$ we write (\ref{CRglobal}) in the form $z_\eta = J(Z) Z_\xi$. This equation can be in turn written as

\begin{eqnarray}
\label{holomorphy}
z_{\bar\zeta} - A(z)\bar z_{\bar\zeta} = 0,\quad
\zeta\in\D.
\end{eqnarray}
Here a smooth map $A: \Omega \to Mat(n,\C)$ is defined by the equality $L(z) v = A \overline v$ for any vector $v \in \C^n$ and $L$ is an $\R$-linear map defined by $L = (J_{st} + J)^{-1}(J_{st} - J)$. It is easy to check that the condition $J^2 = -Id$ iq equivalent to the fact that $L$ is $\overline\C$-linear. The matrix $A(z)$ is called {\it the complex matrix} of $J$ in the local coordinates $z$. Locally the correspondence between $A$ and $J$ is one-to-one. Note that the condition $J(0) = J_{st}$ means that $A(0) = 0$.

If $Z'$ are other local coordinates and $A'$ is the corresponding complex matrix of $J'$, then, as it is easy to check, we have the following transformation rule:

\begin{eqnarray}
\label{CompMat}
A' = ({z'}_zA  + { z'}_{\overline z})({\overline z'}_{\overline z} + {\overline z'}_{ z}A)^{-1}
\end{eqnarray}
(see \cite{SuTu}).

\subsection{Plurisubharmonic functions on almost complex manifolds: the background}  Let  $u$ be a real $C^2$ function on an open subset $\Omega$ of an almost complex manifold  $(M,J)$. Denote by $J^*du$ the
 differential form acting on a vector field $X$ by $J^*du(X):= du(JX)$. Given point $p \in M$ and a tangent vector $V \in T_p(M)$ consider  a smooth vector field $X$ in a
neighborhood of $p$ satisfying $X(p) = V$. 
The value of the {\it complex Hessian} ( or  the  Levi form )   of $u$ with respect to $J$ at $p$ and $V$ is defined by $H(u)(p,V):= -(dJ^* du)_p(X,JX)$.  This definition is independent of
the choice of a vector field $X$. For instance, if $J = J_{st}$ in $\C$, then
$-dJ^*du = \Delta u d\xi \wedge d\eta$; here $\Delta$ denotes the Laplacian. In
particular, $H_{J_{st}}(u)(0,\frac{\partial}{\partial \xi}) = \Delta u(0)$.

Recall some  basic properties of the complex Hessian (see for instance, \cite{DiSu}):

\begin{lemma}
\label{pro1}
Consider  a real function $u$  of class $C^2$ in a neighborhood of a point $p \in M$.
\begin{itemize}
\item[(i)] Let $F: (M',J') \longrightarrow (M, J)$ be a $(J',
  J)$-holomorphic map, $F(p') = p$. For each vector $V' \in T_{p'}(M')$ we have
$H_{J'}(u \circ F)(p',V') = H_{J}(u)(p,dF(p)(V'))$.
\item[(ii)] If $f:\D \longrightarrow M$ is a $J$-complex disc satisfying
  $f(0) = p$, and $df(0)(\frac{\partial}{\partial \xi}) = V \in T_p(M)$ , then $H_J(u)(p,V) = \Delta (u \circ f) (0)$.
\end{itemize}
\end{lemma}
 Property (i) expresses the holomorphic invariance of the complex Hessian. Property (ii) is often useful in order to compute the complex Hessian  on a given 
tangent vector $V$.

 Let $\Omega$ be a domain $M$. An upper semicontinuous function $u: \Omega \to [-\infty,+\infty[$ on $(M,J)$ is
{\it $J$-plurisubharmonic} (psh) if for every $J$-complex disc $f:\D \to \Omega$ the composition $u \circ f$ is a subharmonic function on $\D$. By Proposition
\ref{pro1}, a $C^2$ function $u$ is psh on $\Omega$ if and only if it has    a 
positive semi-definite complex Hessian on $\Omega$ i.e. $H_J(u)(p,V) \geq 0$  for any $ p \in \Omega $ and $V \in T_p(M)$. 
A real $C^2$ function $u:\Omega \to \R$ is called {\it strictly $J$-psh} on $\Omega$, if $H_J(u)(p,V) > 0$ for each $p \in M$ and $V \in T_p(M) \backslash \{ 0\}$. Obviously, these notions   are local: an upper semicontinuous (resp. of class $C^2$) function on $\Omega$ is  $J$-psh (resp. strictly) on $\Omega$ if and only if it is $J$-psh (resp. strictly) in some open neighborhood of each point of $\Omega$.

A useful  observation is that the Levi form of a function $r$ at
a point $p$ in
an almost complex manifold $(M,J)$ coincides with the Levi form with respect
to the standard structure $J_{st}$ of $\R^{2n}$ if {\it suitable} local
coordinates near $p$ are choosen. Let us explain how to construct these adapted
coordinate systems.

As above, choosing local coordinates near $p$ we may identify a neighborhood
of $p$ with a neighborhood of the origin and assume that $J$-holomorphic discs
are solutions of (\ref{holomorphy}).

\begin{lemma}
\label{normalization}
There exists a  second order polynomial local diffeomorphism $\Phi$ fixing the
origin and with linear part equal to the identity such that in the new coordinates
 the complex matrix  $A$  of $J$ (that is $A$ from the equation (\ref{holomorphy})) satisfies
\begin{eqnarray}
\label{norm}
A(0) = 0, A_{z}(0) = 0
\end{eqnarray}
\end{lemma}
Thus, by a suitable local change of coordinates one can remove the 
terms linear in $z$ in the matrix $A$. We stress that in general it is impossible to
get rid of  first order terms containing  $\overline z$ since this would
impose a restriction on the Nijenhuis tensor $J$ at the origin.

I have learned this result  from  unpublished E.Chirka's notes; see .\cite{DiSu} for the proof. In
\cite{SuTu} it is shown that, in an almost complex manifold of (complex) dimension 2,
 a similar normalization is possible along a given embedded $J$-holomorphic disc.
 
 As a typical consequence consider a totally real manifold $E$ defined by (\ref{edge}). Then the function $u = \sum_{j=1}^n \rho^2$ is strictly $J$-psh in a neighborhood of $E$. Indeed, it suffices to choose local coordinates near $p \in M$ according to Lemma \ref{normalization}. This reduces the verification to the well-known case of $J_{st}$.

\section{Filling a wedge by complex discs}  Here we construct a family of Bishop's type discs filling a wedge with a totally real edge. Each disc is glued to the edge $E$ along the upper semi-circle.

\subsection{Case of the standard structure}  First consider the model case $M = \C^n$ with $J = J_{st}$ and $E = i\R^n = \{ x_j = 0, j = 1,...,n \}$. Denote by $W$ the wedge $W = \{ z= x +iy: x_j < 0, j=1,....,n \}$.

Let 
\begin{eqnarray*}
P_0\phi = \frac{1}{2\pi i} \int_{b\D} \phi(\omega) \frac{d\omega}{\omega}
\end{eqnarray*}
denotes the average of a real function $\phi$ over $b\D$, and let also
\begin{eqnarray*}
S \phi(\zeta) = \frac{1}{2\pi i} \int_{b\D} \frac{\omega + \zeta}{\omega - \zeta }\phi(\omega) \frac{d\omega}{\omega}
\end{eqnarray*}
be the Schwarz integral. In terms of the Cauchy transform 

\begin{eqnarray*}
Kf(\zeta) = \frac{1}{2 \pi i} \int_{b\D} \frac{f(\omega)d\omega}{\omega - \zeta}
\end{eqnarray*}
we have the following relation: $S = 2K - P_0$. As a consequence the boundary properties of the Schwarz integral are the same as the classical properties of the Cauchy integral.

For a {\it non-integer} $r > 1$ consider the Banach spaces $C^r(b\D)$ and $C^r(\D)$ (with the usual H\"older norm). This is classical that $K$ and $S$ are bounded linear mappings in these classes of functions. For real function $\phi \in C^r(b\D)$ the Schwarz integral $S\phi$ is a function of class $C^r(\D)$ holomorphic in $\D$; the trace of its real part on the boundary coincides with $\phi$ and its imaginary part vanishes at the origin. In particular,  every holomorphic function $f \in C^r(\D)$ satisfies the Schwarz formula $f = S\Re f + iP_0f$.

We are going to fill $W$  by complex discs glued to $i\R^n$ along the (closed) upper 
semi-circle $b\D^+ = \{ e^{i\theta}: \theta \in [0, \pi] \}$; let also $b\D^-:= b\D \setminus b\D^+$.

Fix a smooth real function $\phi: b\D \to \R$ such that $\phi \vert b\D^+ = 0$ and $\phi \vert b\D^- <  0$.

Consider now  a real $2n$-parameter  family of holomorphic discs $z^0 = (z_1^0,...,z_n^0) : \D \to \C^n$ with components

\begin{eqnarray}
\label{BP1}
z_j^0(c,t)(\zeta) = x_j(\zeta) + iy_j(\zeta) = t_j S\phi(\zeta) + ic_j, j= 1,...,n
\end{eqnarray}
Here $t_j > 0$ and $c_j \in \R$ are  parameters, $t= (t_1,...,t_,)$, $c= (c_1,...,c_n)$. The following properties of this family are obvious:

\begin{itemize}
\item[(i)] for every $j$ one has $x_j \vert b\D^+ = 0$ and $x_j(\zeta) < 0$ when $\zeta \in \D$ (by the maximum principle for harmonic functions). 
\item[(ii)] the evaluation map $Ev_0: (t,c) \mapsto z^0(c,t)(0)$ is one-to-one from $\{ (c,t): c_j > 0, t_j \in \R\}$ on $W$.
\end{itemize}

Our goal is to construct a local  analog of this family in the general case.

\subsection{General case} In order to write an integral equation defining a required family of discs we need to employ an analog of the Schwarz formula and to choose suitable local coordinates . We proceed in several steps.

{\it Step 1.} Recall that for any complex function $f \in C^r(\D)$ the Cauchy-Green transform is defined by

\begin{eqnarray*}
Tf(\zeta) = \frac{1}{2 \pi i} \int\int_{\D} \frac{f(\omega)d\omega \wedge d\overline\omega}{\omega - \zeta}
\end{eqnarray*}
This is classical that $T: C^r(\D) \to C^{r+1}(\D)$ is a bounded linear operator for every non-integer $r > 0$. Furthermore,  $(Tf)_{\overline\zeta} = f$ i.e. $T$ solves the $\overline\partial$-equation in the unit disc. Recall also that the function $Tf$ is holomorphic on $\C \setminus \overline\D$.

 We have the following Green-Schwarz formula (see the proof, for example, in  \cite{SuTu}, although, of course, it can be found in the vaste list of classical works). Let $f = \phi + i\psi : \overline\D \to \C$  be a function of class $C^{r}(\D)$. Then for each $\zeta \in \overline\D$ one has
\begin{eqnarray}
\label{GreenSchwarz}
f(\zeta) = S \phi(\zeta) + i P_0 \psi + Tf_{\overline\zeta}(\zeta) - \overline{Tf_{\overline\zeta}(1/\overline\zeta)}
\end{eqnarray}
Note that because of the "symmetrization" the real part of the sum of two terms containing $T$ vanishes on the unit circle. Notice also that the last term is holomorphic on $\D$.

{\it Step 2.} Now let $(M,J)$ be an almost complex manifold of complex dimension $n$ and $E$ be a totally real $n$-dimensional submanifold of $M$.  We assume that $E$ and $W$ are given by (\ref{edge}) and (\ref{wedge}) respectively.

First, according to Section 2 we choose local coordinates $z$ such that $p = 0$ and the complex matrix $A$ of $J$ satisfies (\ref{norm}). For every $\tau > 0$ small enough and $C > 0$ big enough the functions 
$$\tilde\rho_j:=  \rho_j - \tau \sum_{ k \neq j} \rho_k  + C\sum_{k=1}^n \rho_k^2$$
are strictly $J$-psh in a neighborhood of the origin and the "truncated" wedge $W_\tau = \{ \tilde\rho_j < 0, j=1,...,n \}$ is contained in $W$. After a $\C$-linear (with respect to $J_{st}$) change of coordinates 
one can assume that $\tilde\rho_j = x_j + o(\vert z \vert)$. Consider now a local diffeomorphism 
$$\Phi: z = x_j + iy_j \mapsto z' = x_j' + iy_j' = \tilde\rho_j + iy_j$$
Then $\Phi(0)= 0$, $d\Phi(0) = 0$ and in the new coordinates $\tilde\rho_j = x_j$ (we drop the primes), $E = i\R^n$ and $W_\tau = \{ x_j < 0 , j= 1,...,n\}$. We keep the notation $J$ for the direct image 
$\Phi_*(J)$. Then in the choosen coordinates the complex matrix still satisfies  $A(0) = 0$. Note also that the coordinate functions $x_j$ are strictly psh for $J$. 

Finally, similarly to the proof of Lemma \ref{lemma1}, for $\lambda > 0$  consider the isotropic dilations $d_\lambda: z \mapsto \lambda^{-1}z$ and the direct images $J_\lambda:= (d_\lambda)_*(J)$. Denote by $A(z,\lambda)$ the complex matrix of $J_\lambda$.

{\it Step 3.} For $\lambda > 0$ small enough we are looking for the solutions $z: \D \to \C^n$ of Bishop's type integral equation

\begin{eqnarray}
\label{MainBishop}
z(\zeta) = h(z(\zeta), t,c,\lambda) 
\end{eqnarray}
with
$$h(z(\zeta), t,c,\lambda) = tS \phi(\zeta) + i c + TA(z,\lambda) z_{\zeta}(\zeta) - \overline{T A(z,\lambda) z_{\zeta}(1/\overline\zeta)}$$
where $t = (t_1,...,t_n)$, $t_j > 0$ and $c \in \R^n$ are real parameters. Note that the first and the last terms in the right hand are holomorphic on $\D$. Therefore, any solution of (\ref{MainBishop}) satisfies the Cauchy-Riemann equations (\ref{holomorphy}) i.e. is a $J$-complex disc. Furthemore, $x_j(\zeta)$ vanishes on $b\D^+$ ( that is $z(b\D^+) \subset E$) and is negative on $b\D^-$. Since the function $z \mapsto x_j$ is strictly $J$-psh, by the maximum principle the image $z(\overline\D)$ is contained in $\overline W_\tau$.

The existence of solutions follows by the implicit function theorem. Note that for $\lambda = 0$ the equation ( \ref{MainBishop}) admits the solution ( \ref{BP1}). Consider the smooth map
of Banach spaces

$$H: C^r(\D) \times \R^n \times \R^n \times \R  \longrightarrow  C^r( \D)$$
$$ H: (z, c,t, \lambda)  \mapsto h(z(\zeta), t,c,\lambda)$$ 

Obviously the partial derivative of $H$ in $z$ vanishes: $(D_zH)(z^0,c,t,0) = 0$, where $z^0$ is a disc given by ( \ref{BP1}). By the implicit function theorem, for every $(c,t, \lambda)$ close enough to the origin the equation ( \ref{MainBishop}) admits a unique solution $z(c,t)( \zeta)$ of class $C^r( \D)$, smoothly depending on parameters $(c,t)$  (as well as  $\lambda$, of course). 

Fixing $ \lambda > 0$, consider the smooth evaluation map
$$Ev_ \lambda: (c,t)  \mapsto z(c,t)(0)$$
which associates to each parameter the centre of the corresponding disc. Note that for $ \lambda = 0$ we obtained the linear mapping $Ev_0(c,t) = (tS \phi(0) + ic)$ appeared already  in the model case of the standard structure. More precisely, for every $\lambda > 0$ and  every $c \in \R^n$     the equation (\ref{MainBishop}) has the unique constant solution $z = ic$ when $t = 0$. Hence, $Ev_\lambda$ admits the expansion $Ev_\lambda(c,t) = Ev_0(c,t) + o(\vert (c,t,\lambda)\vert)$.

Denote by $V$ the wedge $V = \{ (c,t): c \in \R^n, t_j > 0, j=1,...,n \}$ with the edge $L = \{ (c,0): c \in \R^n \}$. Then $Ev_0(V)$ coincides with $W_\tau$ in a neighborhood of the origin.  Furthermore, 
$Ev_\lambda(L) = E$ for $\lambda > 0$. Also, for $\alpha > 0$ the "truncated" wedge 
$W_\alpha = \{ x_j - \alpha\sum_{k \neq j} x_k < 0 \}$ with the edge $E$ is contained in $W_\tau$. The faces of the boundary of $W_\alpha$ are transverse to the face of $W_\tau$. Since this property is stable under small perturbations,  we conclude that $W_\alpha \subset Ev_\lambda(V)$ for all $\lambda$ small enough. In terms of the initial defining functions $\rho_j$ we have $\{ z: \rho_j - \delta\sum_{k \neq j} \rho_k < 0 \} \subset W_\varepsilon$ when $\tau + \varepsilon < \delta$.

We sum up in the following

\begin{thm}
\label{theo2}
Let $E$ be a totally real submanifold (\ref{edge}) of an almost complex manifold $(M,J)$ and $W$ be a wedge (\ref{wedge}). For every $\delta > 0$ there exists a family of $J$-complex discs smoothly depending on $2n$ real parameters with the following properties:
\begin{itemize}
\item[(i)] the boundary of every disc is glued along $b\D^+$ to $E$ and every disc is contained in $W$;
\item[(ii)] every point of the truncated wedge $\{ z: \rho_j - \delta\sum_{k \neq j} \rho_k < 0 \}$ belongs to some disc.
\end{itemize}
\end{thm}
Although we do not need it in this paper,  note that every point of the wedge $E$ belongs to the boundary of some disc. Concerning regularity of manifolds and almost complex structures, it suffices to require the class $C^\alpha$ with real $\alpha > 2$. We skip the details.

\section{Boundary uniqueness}

Now we can easily prove the main result of the present paper.

\begin{thm}
\label{theo1}
Let $E$ be a totally real submanifold (\ref{edge}) of an almost complex manifold $(M,J)$ and $W$ be a wedge (\ref{wedge}).  Assume that $u$ is a plurisubharmonic function in  $W$ such that 
\begin{eqnarray}
\label{boundary1}
\lim \sup_{W \ni z \to E} u(z) = -\infty
\end{eqnarray}
Then $u \equiv -\infty$.
\end{thm}
Indeed, for every disc $f$ constructed in Theorem \ref{theo2} the function  $(u \circ f)(\zeta)$  is subharmonic on $\D$ and tends  to $-\infty$ when $\zeta \to b\D^+$; therefore, $u \circ f \equiv -\infty$ 
(by classical properties of subharmonic functions, see \cite{R}). Now by (ii) Theorem \ref{theo2} we conclude that $u \equiv -\infty$.

\begin{cor}
Let $E$ be a closed subset of an almost complex $n$-dimensional manifold $(M,J)$. Suppose 
that $E$ contains a germ of a totally real manifold of dimension $n$. Then $E$ is not a pluripolar set.
\end{cor}

\end{document}